\documentclass[12pt]{amsart}
\usepackage {amsmath, amscd, amsbsy, amsfonts}
\usepackage{amsmath,latexsym,amssymb,amsfonts,mathrsfs,bbm}
\usepackage{fullpage}
\usepackage{graphicx}
\usepackage{hyperref}
\newlength{\baseunit}               
\newcommand{\R}{\mathbbmss{R}}

\newcommand{\s}{\mathbf{s}}

\def \gam {{\Gamma}}

\def \n {{\mathcal N}}
\def \M {\mathcal M}
\def \D {\mathcal D}
\def \Q {\mathbb Q}
\def \proj {{\mathbb{P}^r}}

\def \d {\Delta}

\def \O {{\mathcal O}}
\def \H {{\mathcal H}}

\def \G {{\mathcal G}}

\def \F {{\mathcal F}}

\def \s {{\mathcal S}}
\def \R {{\mathcal R}}
\def \C {{\mathcal C}}
\def \Cu {{\mathcal C \mathcal U}}
\def \T {{\mathcal T}}
\def \n {{\mathcal N}}
\def \mbar {\overline{\mathcal M}}

\def \wt {\widetilde}
\def \E {{\mathcal E}}

\def \J {{\mathcal J}}

\def \gone {\mbar_{1,n}(r,d)^*}
\def \gof {\mbar_{1,4}(r,d)^*}
\def \mof {\overline{M}_{1,4}}
\def \mon {\overline{M}_{1,n}}
\def \mb {\overline{M}}
\def \pba {\mathbb P^3}
\def \pbon {\mathbb P^4}
\def \pnam {\mathbb P^5}
\def \pr {\mathbb P^r}
\def \mmnrd {\mbar_{1,n}(r,d)}
\def \beq {\begin{eqnarray*}}
\def \eeq {\end{eqnarray*}}            
\def \one {\mathbf 1}    
\def \ss {\mathbb S}

\newtheorem{theorem}{Theorem}[section] 

\newtheorem{lemma}[theorem]{Lemma}
\newtheorem{example}[theorem]{Example}

\newtheorem{proposition}[theorem]{Proposition}

\newcommand{\bpf}{\noindent {\it Proof.  }}
\newcommand{\epf}{\qed \vspace{+10pt}}

\title[Enumerative Invariants]{Characteristic numbers of elliptic space curves}

\author{Dung Nguyen}

\begin{document}
\maketitle
\begin{abstract}
We solve the problem of characteristic
numbers of elliptic curves in any dimensional
projective space. The answers are given
in the form of effective recursions. A C++ program
implementing most of the recursions 
is available upon request. 
 \vspace{-20pt}
\end{abstract} 
\section{introduction}
 Computing the charateristic numbers of curves in projective spaces is a classical problem in algebraic
geometry: how many curves in $\proj$ of given degree and genus 
that pass through a general set
of linear subspaces, and are tangent to a general set of hyperplanes? 
These types of problems have provided a wealth of inspiration
for classical algebraic geometers, such as Schubert, Zenthen,
Mailard, and have been one of the main drives of the 
development of intersection theory.  They have achieved remarkable
success despite the technology of their time (in mathematical
and literal sense). For example, Schubert has correctly
computed the number of twisted cubics in $\mathbb P^3$ that
are tangent to $12$ quadric surfaces ($5,819,539,783,680$), and
this represents the summit of enumerative geometry at the time (see \cite{ksx}).

 However, not until recently are characteristic numbers computed
in certain generality, as all of the results obtained by classical 
algebraic geometers are bounded by degree. Moduli space of stable maps
provide a powerful tools to attack characteristic numbers problems
of curves with low genus. The case of rational space curve of
any degree was solved
in \cite{idq}. The author actually gave algorithms to compute
intersections of divisors on $\mbar_{0,n}(r,d)$ and derived
characteristic numbers as a corollary. However, if one is only
interested in the characteristic numbers, then a small tweak
would give a much faster and easier to implement algorithm. The
case of elliptic plane curves, any degree, was solve in
\cite{char}. The difficulty to extend that result to
$\mathbb P^r$ was due to two factors. Firstly, there was yet
a result calculating incidence-only (no tangency) numbers for elliptic
space curves in any dimensional projective space. Secondly, because
the tangency divisor on $\mbar_{1,n}(r,d)$ is related
to the divisor of elliptic curves with fixed $j$-invariant,
hence one also needs to compute the characteristic numbers
of elliptic curves with fixed $j$-invariant. The latter
was solved in a separate paper by the author (\cite{dn}). The
algorithms for characteristic numbers of genus two
plane curves were given in \cite{gkp}, although no actual 
numbers were shown.

 In this paper, we solve the full characteristic numbers problem for
elliptic curves in projective spaces by using the moduli space of
stable maps. Note
that, even in the incidence-only case our result is already new. For example, no
incidence-only characteristic numbers for elliptic curves
in $\mathbb P^4$ and $\mathbb P^5$ have been computed. The number
of elliptic curves in $\mathbb P^3$ were computed in \cite{eg}
and the proof was scheduled to appear in another paper. However, the author
of this paper was unable to locate it. The incidence-only numbers
of elliptic curves in $\mathbb P^3$ were also computed in
\cite{ratell}. In that paper, only the numbers upto degree $4$ elliptic
curves were computed, but the method could compute
numbers of any degree. The incidence-only numbers of degree $4$ elliptic
curves in $\mathbb P^3$ were first computed in \cite{av} using
a classical argument.
 
 The approach of this paper is as follows. We first compute
incidence-only numbers for elliptic space curves using
Getzler relation. This is analogous to the enumeration
of rational space curves using WDVV equation, but with a
small twist. First we will not be using the entire
moduli space of stable maps $\mmnrd$ as intersecting
enumerative classes with the virtual fundamental class
will have unwanted contributions from other components (this
space is not irreducible). In fact, it was claimed in
\cite{eg} that the actual count of elliptic curves is
a linear combination of genus one and genus zero Gromov-Witten invariants.
We will only use the main component, that is the closure
of the locus of maps with smooth source curves, denoted by
$\mbar_{1,r}(r,d)^*$. For our purpose,
knowledge of the maps on the loci where the main component
intersect other components is required, and this is desribed
in Section $4$.

In the case of rational curve, the WDVV equation
is a rational equivalence of divisors, hence we can pull back
via the forgetful morphism
and obtain a relation on the stable map space. This is no longer
true for the genus one stable map space. The reason is that Getzler relation
is a rational equivalence of codimension $2$ strata on $\mof$, and that the 
forgetful morphism $\mbar_{1,4}(r,d)^* \to \overline{\M}_{1,4}$
have fibre dimensions that could jump : for example, the preimage
of the stratum $\delta_{2,2}$, has two components, one of which is
of codimension $2$, and the other is a Weil divisor. We can get
around this as follows. First we use enough enumerative constraints
to cut down the space $\gone$ into a $2$-dimensional family. Then
we pushforward via the forgetful morphism, and then intersect with
 Getzler's relation. As a result we obtain a relation of the enumerative
invariants of elliptic curves with that of rational curves,
rational cuspidal curves, and elliptic curves with fixed $j$ invariants.
The first is well-known, and the latter two were computed in
\cite{dn2} and \cite{dn} respectively. 

 To go from incidence-only to full characteristic numbers,
we study the relation between the tangency divisor and
the incident and boundary divisors. We will not
obtain a rational equivalence, due to the present
of enumeratively irrelevant divisors, those that 
are intersections of $\gone$ with other components. However,
we obtain a numerical equivalence whenever we intersect with
only curves in $\gone$ that has empty intersection with the
irrelevant divisors. This is true for $1$-dimensional
families of elliptic curves cut down by enumerative constraints,
so this is enough for our purpose.

 The structure of this paper is as follows. In Section $2$, we introduce
basic notions such as various stacks of stable maps and the enumerative
constraints. In Section $3$ we review Getzler's relation
and give an example of using the relation to count elliptic curves.
In Section $4$, we give the recusion computing incidence-only
numbers of elliptic space curves based on Getzler's relation.
In Section $5$, we give the recursions computing full
characteristic numbers of elliptic space curves. We end with
some tables with numerical examples in Section $6$.

 The author would like to thank Ravi Vakil for many helpful suggestions 
and some of the arguments used in this paper, and for suggisting this
problem.

\section{Definitions and Notations}
\subsection{The moduli space of rational and elliptic curves in $\proj$.}
We denote $\mbar_{0,n}(r,d)$ the Kontsevich compactification of the moduli
space of genus zero curves with $n$ marked points of degree $d$ in $\proj$. Let
$\gone$ be the main component of the moduli space of stable maps
of genus one, that is, the closure of the locus of maps with smooth
domains.
 We will also use the notation $\mbar_{0,S}(r,d)$
and $\mbar_{1,S}(r,d)$ where the markings are indexed
by a set $S$.

\subsection{The constraints and the ordering of constraints.}

 We will be concerned with the number of curves passing through a constraint, and each 
constraint is denoted by a $(r+1)$-tuple $\d$ as follows : \\
{\bf (i)}  $\d(0)$ is the number of hyperplanes that the curves need to be tangent to. \\
{\bf (ii)} For $0<i \leq r$, $\d(i)$ is the number of subspaces of codimension
$i$ that the curves need to pass through. 

 In \cite{dn} and \cite{dn2}, the constraints may have $(r+2)$ coordinates
because we want to impose conditions on the node or cusp, but we will not
need that here.

 Note that because in general a curve of degree
$d$ will always intersect a hyperplane at $d$ points, introducing an incident condition
with a hyperplane has the same effect as that of multiplying the enumerative number by $d$. 
For example, if we ask how many genus zero curves of degree $4$ in $\mathbb P^3$ that pass through the constraint
$\d = (1,2,3,4,0)$, that means we ask how many genus zero curves of degree $4$ pass through
three lines, four points, are tangent to one hyperplane, and then multiply that answer
by $4^2$. We will also refer to $\d$ as a set of linear spaces,
hence we can say, pick a space $a$ in $\d$. 

 We consider the following ordering on the set of constraints, in order to prove that our algorithm will 
terminate later on. Let $r(\d) = -\sum_{i >1}^{i\leq r} \d[i]\cdot i^2$, and this will be our rank function.
We compare two constraints $\d,\d'$ using the 
following criteria, whose priority are in the following order. We only proceed
to using the next criterion if using the current
one give us a tie.
\begin{itemize}
\item If $\d(0) > \d'(0)$  then $\d < \d'$.
 \item  If $\d(0) = \d'(0)$ and $\d$ has fewer non-hyperplane elements than $\d'$ does, then $\d<\d'$. 
\item If $r(\d) < r(\d')$ then $\d < \d'$. 
\end{itemize}

Informally speaking, characteristic numbers
where the constraints are more spread out at two ends 
are computed first in the recursion.
We write $\d = \d_1\d_2$ if $\d = \d_1 + \d_2$
as a parition of the set of linear spaces in $\d.$

\subsection{The stacks $\R,\n, \C, \n\R, \R\R, \R\R_2, \E, \E\R, \E\R\R, \R\E\R$ .}
 We list the following definitions of stacks of stable maps that will
occur in our recursions. \\ \\
{\bf 1)} Let $\R(r,d)$  be the usual moduli space of genus zero stable maps, 
$\mbar_{0,0}(r,d)$. \\ \\
{\bf 2)} Let $\n(r,d)$ be the closure in $\mbar_{0,\{A,B\} }(r,d)$ of the locus of maps of smooth
rational curves $\gamma$ such that $\gamma(A) = \gamma(B)$. Informally,
$\n(r,d)$ parametrizes degree $d$ rational nodal curves in $\proj$. \\ \\
{\bf 3)} Let $\Cu(r,d)$ be the closure in $\mbar_{0,\{C\}}(r,d)$ of
the locus of maps of smooth rational curves $\gamma$ such that 
the differential $\gamma'(C)$ is zero. Informally, $\C(r,d)$
parametrizes degree $d$ rational cuspidal curves in $\proj.$\\ \\
{\bf 4)} Let $\E(r,d)$ be the main component of the moduli space
of genus one stable maps $\mbar_{1,0}(r,d)$. That is
$\E(r,d) = \mbar_{1,0}(r,d)^*$. \\ \\
{\bf 5)} For $d_1,d_2>0,$ let $\R\R(r,d_1,d_2)$ be $\mbar_{0,\{C\}}(r,d_1) \times \mbar_{0,\{C\}}(r,d_2)$
where the fibre product is taken over evaluation maps $ev_{C}$ to $\proj.$ \\ \\
{\bf 6)} Similarly we can define $\n\R(r,d_1,d_2)$, $\E\R(r,d_1,d_2)$ (see figure 1).\\ \\
{\bf 7)} For $d_1,d_2>0$, let $\R\R_2(r,d_1,d_2)$ be the closure in $\mbar_{0,\{A,C\}}(r,d_1) \times_{\proj} \mbar_{0,\{B,C\}}(r,d_2)$  
(the projections are evaluation maps $e_C$) of the locus of maps $\gamma$ 
such that $\gamma(A) = \gamma(B)$. We call maps in this family
rational two-nodal reducible curves.
\\ \\ 
{\bf 8)} For $d_1,d_2,d_3>0$, let $\R\E\R(r,d_1,d_2,d_3)$ be $\mbar_{0,\{C\}}(r,d_1) \times_{ev_C} \mbar_{\{1,C,D\}}(r,d_2) \times_{ev_D} \mbar_{0,D}(r,d_3) $.
Similarly, we can define $\E\R\R(d_1,d_2,d_3)$ for $d_1,d_2,d_3 > 0$. \\ \\
{\bf 9)} We define $\J(r,d)$ to be the closure in $\E(r,d)$ of the locus of maps
whose domains are smooth and have a fixed but generic $j$-invariant. The enumerative
geometry of this stack is studied in \cite{dn}. \\ \\ \\
$$\includegraphics[width = 150mm]{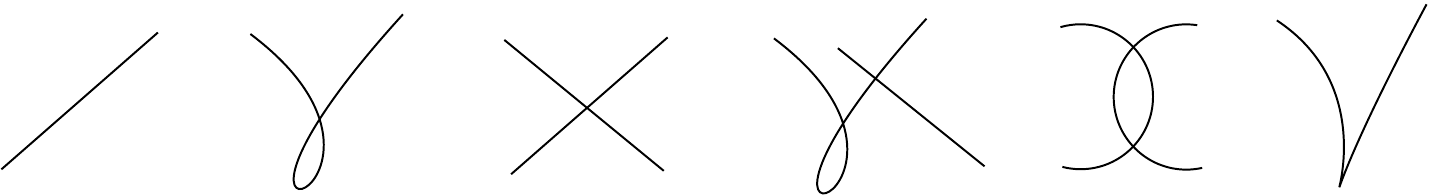}$$ \\
$$\text{Fig 1. Pictorial description of a general curve in the stacks $\R, \Cu, \n, \R\R, \n\R, \R\R_2.$}$$
$$ $$
$$ $$
$$ $$
$$ $$
$$ $$
$$\includegraphics[width = 150mm]{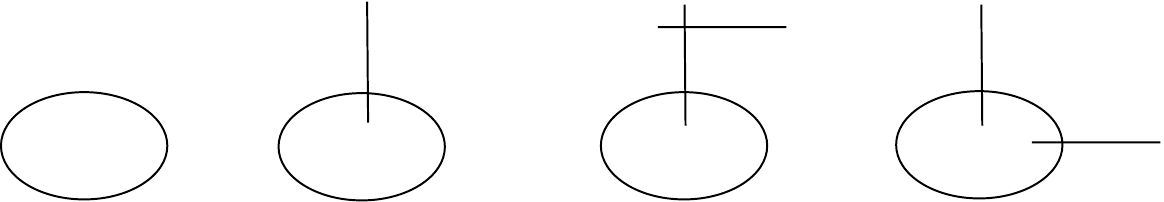}$$ \\
$$\text{Fig 2. Pictorial description of a general curve in the stacks $\E, \E\R, \E\R\R, \R\E\R.$}$$

\subsection{Stacks of stable maps with constraints.} Let $\F$
be a stack of stable maps of curves into $\proj$. For a constraint
$\d$, we define
$(\F, \d)$ be the closure in $\F $ of the locus of maps
that satisfy the constraint $\d$. 
For maps of reducible source curves, tangency
condition include the case where the image of the node
lies on the tangency hyperplane, as the intersection
multiplicity is $2$ in this case. For a stack $\F$ that
is supported on a finite
number of points then we denote
$\# \F$ to be the stack-theoretic length of $\F$. 

If $\F$ is a closed substack of the stacks $\n\R, \R\R, \E\R$ then we denote
$(\F,\gam_1,\gam_2,k)$  to be the closure in $\F$
of the locus of maps $\gamma$ such that the
restriction of $\gamma$ on the $i-$th component
satisfies constraint $\gam_i$, and that
the image of the node lies on a codimension
$k$ subspace. We use the notation $(\F,\d,k)$ if we don't want
to distinguish the conditions on each component. If $k$ is $0$
we omit it from the notation.

If $\F$ is a closed substack of the stack $\R\R_2$ then we denote
$(\F,\gam_1,\gam_2)$  to be the closure in $\F$
of the locus of maps $\gamma$ such that the
restriction of $\gamma$ on the $i-$th component
satisfies constraint $\gam_i$. We use the notation $(\F,\d)$.

If $\F$ is a closed substack of the stacks $\E\R\R,\R\E\R$ 
then we denote
$(\F,\gam_1,\gam_2, \gam_3)$  to be the closure in $\F$
of the locus of maps $\gamma$ such that the
restriction of $\gamma$ on the $i-$th component
satisfies constraint $\gam_i$. We use the notation $(\F,\d)$ if we don't want
to distinguish the conditions on each component.

 The enumerative geometry of all the stacks defined
above are known, except
for the stacks that involve $\E$. But the enumerative
geometry of the stacks $\E\R,\R\E\R,\E\R\R$ can be
easily deduced from that of $\E$ (see \cite{dn}, Section
$3$). Note the small difference between the notation
here and in \cite{dn}, \cite{dn2}, as we have remove
some conditions on the nodes as they are not necessary.

\section{Getzler's relation on $\mof$}
  We review some of the basic intersection
theory on $\mof$ and especially Getzler's
relation. For a more complete treatment
see \cite{eg}. Consider the moduli
space $\mof$ with the $\ss_4$ action
by permuting the marked points. In \cite{eg},
the following $\ss_4$-invariant codimension
$2$ strata are defined:

$$\includegraphics[width = 130mm]{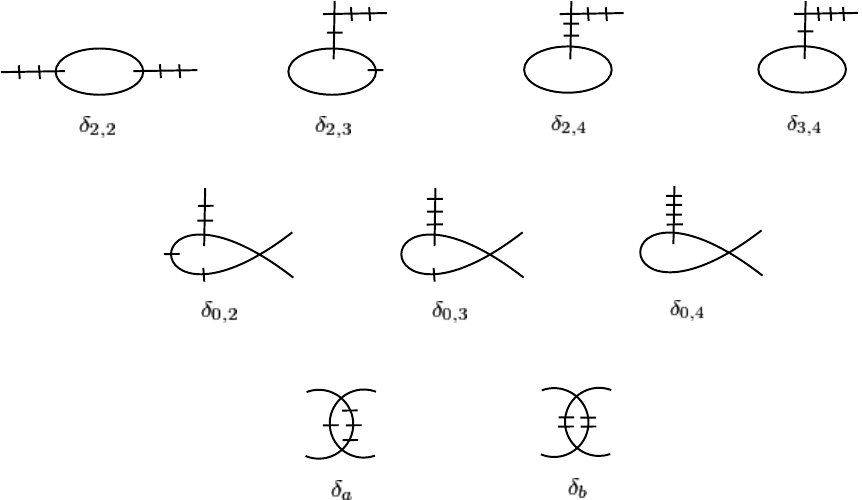}$$
$$\text{Fig 3. The $S_4$-invariant codimension $2$ strata of $\mof$}$$

 Getzler computed the intersection matrix of
these cycles and found a non-trivial null vector
$$12\delta_{2,2} - 4\delta_{2,3}-2\delta_{2,4}+6\delta_{3,4} + \delta_{0,3}+\delta_{0,4}-2\delta_{\beta} \cong 0 $$

These cycles are understood stack-theoretically. Pandharipande gave a beautiful argument in \cite{pg} showing
that the above relation is actually a rational equivalence, by
using the auxilliary space of admissible covers.

 Note that it is easy to get relations on $\mon$ induced
by WDVV. Let
$\delta_0$ be the boundary stratum of irreducible
nodal curves on $\mon$. We have a map
$\eta: \mb_{0,n+2} \to \mon$ by identifying the last
two marked points. Pushing forward the WDVV
relation on $\mb_{0,n+2} $ yields a relation on $\mon$. Unfortunately,
these type of relations does not help us counting
elliptic curves cause the general curves
in the strata involed in the relations do
not have a smooth elliptic component.

 Now as an example, we use Getzler's relation to
compute a trivial example.
\begin{example}
 There is $1$ elliptic cubic in $\mathbb P^3$ that passes through
$3$ points and $6$ lines.
\end{example}

\bpf The $3$ point conditions determine the plane of the elliptic
cubic. The other $6$ conditions translate to $6$ other point
conditions for a smooth plane cubic. There is $1$ elliptic
plane cubic through $9$ points is a fact followed easily
from linear algebra. We give another proof using 
Getzler's relation. Call the number of elliptic cubics
satisfying the given constraint $X$. First we consider a $2$-dimensional
of $\F$ maps $\gamma$ on $\mbar_{1,A,B,C,D}(3,3)^*$ satisfying the following
conditions : \\
- $\gamma(A), \gamma(B), \gamma(C), \gamma(D)$ each belongs to a general plane.  \\
- The image of $\gamma$ passes through $3$ point and $4$ lines.
By elementary dimension counting it is easy to see that 
$\F$ is $2$ dimension. Let $p: \gof \to \mof$ be the forgetful
morphism. Now we intersection $\G = p_*(\F)$ with the Getzler's
relation.  \\ \\
- $\G \cap \delta_{2,2}$: There are no map of elliptic
curves with degree $1$. Maps of degree $2$ are double covers,
but with incidence-only constraints, these will not contribute.
Thus the degree on the elliptic components are $3$ and the
two rational twigs are contracted. (We will see in Section
$4$ and $5$ that the maps that contracted the elliptic component
do not contribute either). For each possible way of distributing 
the marked points on the two rational twigs we see that the count is $X$, thus
we have $3X$ in total. Multiply this by the coefficient
of $\delta_{2,2}$ we get $36X.$ \\  \\
- $\G \cap \delta_{2,3}$: Argue similarly, the two 
rational twigs must contract. Then the node on the elliptic
curve must satisfy incidence conditions of $3$ of the points
of $A,B,C,D$ (whichever on the contracted components). Thus we get
the numbers of elliptic curves passing through $4$ points and $4$ lines,
multiplied by $3$ due to a marked point on the elliptic component 
satisfying a hyperplane condition. But this is $0$ cause
an elliptic cubic is planar and cannot pass through
$4$ general points. \\ \\
- Similarly $\G \cap \delta_{2,2}$ and $\G \cap \delta_{2,4}$
are $0$. \\ \\
- $\G \cap \delta_{0,3}$: If the rational twig is contracted,
then we get the numbers of nodal plane cubics passing through
$4$ points and $3$ lines, which is $0$. If the genus $1$ component is contracted,
then by a result in section $4$, the resulting image curve
has a cusp (corresponding to the node on the domain). The cusp
must satisfy an hyperplane incidence condition. By dimension counting,
the contribution in this case is also $0$. \\ \\
- $\G \cap \delta_{0,4}$: Similarly, if the rational twig get
contracted, then there is no contribution. If the genus
$1$ component gets contracted, then we the number
of rational cuspidal cubic curves passing through $3$ points
and $4$ lines. Translating this to the plane curve counting
problem we get the number of rational cuspidal cubic passing
through $7$ points, which is $24$. There are $4$ marked points
each satisfying a hyperplane condition, and each map
has automorphism group of order $2$. Thus the total contribution
is $3^4\cdot 24/2 = 972.$ \\ \\
- $\G \cap \delta_b$: In general, counting maps with
this type of domains ($\R\R_2$) is not trivial, and an algorithm
for this is given in \cite{dn}. In this particular example, however,
it is easy. If any of the twig gets contracted, we have
rational nodal curve with conditions on the node, but
by simple dimension counting we see that there is no contribution.
 Thus the pair has to be line-conic
intersecting at $2$ points. The $3$ points conditions determine
the plane for this pair, so we can translate this into a plane curve 
counting problem. How many pair of line-conics passing through
$7$ points. The answer is $1$, multiplied by the number of ways
to distribute the conditions, which is $\binom{7}{2} = 21$.
 We also have incidence on the marked points,
especially on the conics (each contributes a factor of $2$). For each
marked points distribution, the contribution is $2^2 = 4$. There are
six total ways to distribute the marked points. Thus the total contribution of this 
stratum is $21 \cdot 4 \cdot 6 = 504.$ Multiplying with the coefficient
of the stratum, we get $504 \cdot ( -2) = -1008$. \\ \\
Thus we end up with the equation
$$36X + 972 - 1008 = 0$$
hence $X=1$.
\section{Incidence-only characteristic numbers of elliptic space curves via Getzler's relation}
 First we need a lemma to tell us which maps can occur
as boundary of the main component $\gone$.

The Euler characeristic is an important invariant
in a family of curves. It is $1$ minus the arithmetic
genus. The Euler characteristic can be defined for
any scheme $X$ as the Euler characteristic of its structure
sheaf $\chi(X)$. The following lemma is useful to compute
Euler characteristic of reducible curves.
\begin{lemma}
  Euler characteristic satisfies the inclusion-exclusion
property. That is, if $S_1,S_2$ are two subschemes of an ambient scheme, then
$\chi(S_1 \cup S_2) + \chi(S_1 \cap S_2) = \chi(S_1) + \chi(S_2)$
\end{lemma}
\bpf We have the following exact sequence of sheaf on $S_1 \cup S_2$ :
$$0 \to \O_{S_1\cup S_2} \to \O_{S_1} \oplus \O_{S_2} \to \O_{S_1 \cap S_2} \to 0$$
where all the sheaves on smaller schemes are pushforwarded to $S_1 \cup S_2$.
The identity follows from the fact that Euler characteristic is additive
along an exact sequence. \epf
\begin{proposition}
  Let $\gamma : (C,S) \to \pr$ be a map in
$\mbar_{1,S}(r,d)^*$ ($d\neq 1$) from a
nodal, $n$-pointed curve $C$ of arithmetic genus one 
that contracts a connected union of irreducible 
components having total arithmetic genus one. Let
$R_j$'s ( $1 \leq j \leq m)$ 
be the non-contracting rational tails
and let $n_j$ be the intersection of $R_j$
with the contracted subcurve. Then the
images of the tangent vectors to $R_j$'s
at $n_j$'s are not independent.
\end{proposition}
\bpf First note that contracted subcurve
has arithmetic genus one, hence
each $R_j$ intersects it
at exactly one point. Let $\C \to B$
be a one dimensionsinal family with
general fibers smooth of genus $1$ and
the central fiber being $C$. If $\gamma$
contracts the central fiber then it factors
through a contraction $c : \C \to C'$.
The central fiber of $\C'$, denoted by
$C'$ consists $m$ rational
twigs connected at one node. However,
the tangent vectors to these twigs
at the node are not independent for the
following reason. The Euler characteristic
of the general fiber is $0$ (equals $1-$
the arithmetic genus), and this is a
constant of the fibers in the family. If
the tangent vectors were independent, then
the node (with the intersection scheme
structure ) is reduced, and by
a repeat application of Lemma $4.1 $,
we see that the Euler characteristic
would be $1$. It follows that
the images of the tangent vectors
in $\pr$ are also not independent. \epf

 Using the lemma, one can see that the following
loci are intersections of the main component with
other components. 
For $1 \leq m \leq r, \mbar^m_{1,n}(r,d)$ is the closure of the locus
of maps $\gamma$ having the following property: \\
- The domain of $\gamma$ has a smooth elliptic
component joined with $m$ rational tails $R_j$
at the node $n_j$. \\
- $\gamma$ has degree $0$ on the elliptic component,
and the images of the tangent vectors to $R_j$
at $n_j$ are linearly dependent.

 For $m = 1$, the image curve has a cusp, for $m=2$,
the image curve is a two-component rational curves
meeting at a tacnode. In general , for $m \leq r$ the image curve
has a $m$-fold elliptical singularity. For $m > r$, 
$\mbar^m_{1,n}(r,d)$ lies entirely on $\gone$.

 We are now ready to derive a recursion counting
elliptic space curves. First we need some notations.
  Let $\d$ be a constraint such that $\d(0) = 0$, and assume
$i\geq2$ is an index such that $\d(i) \geq 2$. Let $p',q'$ be
two subspaces in $\d$ of codimension $i$. Let $h,k$ be
two general hyperplanes in $\pr$ and let $p,q$ be subspaces
of codimension $i-1$.
The following constraints
are derived from $\d$ : \\
- $\wt{\d}$ by removing $p'$ and $q'$. \\
- $\d_0$ by removing $p',q'$ but adding a  codimension $i-1$
and a codimension $i+1$ subspaces. \\
- $\d_1$ by removing $p',q'$ but adding a codimension $2$
and a codimension $2(i-1)$ subspaces. \\
- $\d_2$ by removing $p',q'$ but adding a hyperplane and
a codimension $2i-1$ subspace. \\
- $\d_c$ by removing $p',q'$ but adding $p,q,h,k$.
 We denote $\one_{ST}$ be the indicator function where $ST$ is
a logical statement. Let $\ss_4$ acts on the set $\{p,q,h,k\}.$
It is easy to see that $\d_i, 0\leq i \leq 2$ is of lower rank
than that of $\d$.
\begin{theorem}
We have the following recursive formula, providing
the left-hand side is finite (it is understood 
that if a constraint contains the empty subspace then the corresponding
enumerative term is automatically $0$).

\begin{eqnarray*}
 \# (\E(d,r),\d) &=&  \frac{1}{12 (2 + \one_{(i = 2)})} \big [ (-12 \# (\E(r,d),\d_1) \\ 
 &-& \frac{12}{4}  \sum_{d_1 + d_e = d} \sum_{\beta \in \ss_4} \#((\E\R(d_e,d_1,r),\wt{\d}), (\beta(h) \cap \beta(k)), (\beta(p), \beta(q)) \\
  &-& \frac{12}{8} \sum_{d_1+d_2+d_e = d} \sum_{\beta \in \ss_4} \#((\R\E\R(d_1,d_e,d_r),\wt{\d}),(\beta(h),\beta(k)),\emptyset, (\beta(p),\beta(q))) \\
  &+& \frac{4}{2} \sum_{d_1+d_e = d} \sum_{\beta \in \ss_4} \#((\E\R(d_e,d_1,r),\wt{\d}), (\beta(h)), (\beta(k), \beta(p) \cap \beta(q))) \\
  &+& \frac{4}{2} \sum_{d_1+d_e = d} \sum_{\beta \in \ss_4} \#((\E\R(d_e,d_1,r),\wt{\d}), (\beta(h)), (\beta(p),\beta(q)),\beta(k)) \\
  &+& 24 \#(\E(d,r),\d_0) + 24 \#(\E(d,r),\d_2)  \\ 
  &+& \frac{4}{2} \sum_{d_1+d_2+d_e=d} \sum_{\beta \in \ss_4} \#((\E\R\R(d_e,d_1,d_2),\wt{\d}), (\beta(h)), (\beta(k)), (\beta(p),\beta(q))) \\
  &+& \frac{2}{4} \sum_{d_1+d_e = d} \sum_{\beta \in \ss_4} \#((\E\R(d_e,d_1),\wt{\d}),\emptyset, (\beta(h),\beta(k),\beta(p) \cap \beta(q))) \\
  &+& \frac{2}{4} \sum_{d_1 +d_2 + d_e= d} \sum_{\beta \in \ss_4} \#((\E\R\R(d_e,d_1,d_2),\wt{\d}),\emptyset, (\beta(h),\beta(k)), (\beta(p),\beta(q))) \\
  &-& \frac{6}{6} \sum_{d_1+d_e =d} \sum_{\beta \in \ss_4} \#((\E\R(d_e,d_1),\wt{\d}).\emptyset, (\beta(k), \beta(p), \beta(q)),\beta(h)) \\
  &-& \frac{6}{6} \sum_{d_1+d_2+d_e=d} \sum_{\beta \in \ss_4} \#((\E\R\R(d_e,d_1,d_2),\wt{\d}),\emptyset, (\beta(h), (\beta(k), \beta(p), \beta(q)))) \\
  &-& \frac{6}{6} \sum_{d_1+d_n =0} \sum_{\beta \in \ss_4} \frac{1}{2}\#((\n\R(d_n,d_1),\wt{\d}),(\beta(h)),(\beta(k),\beta(p),\beta(q))) \\
  &-&  \sum_{d_1+d_n =d} \frac{1}{2}\#((\n\R(d_n,d_1,r),\wt{\d}),\emptyset, (h,k,p,q))  - \frac{1}{2} \#(\s,\d_c)\\ 
  &+& \frac{2}{8}  \sum_{d_1+d_2 =d}\sum_{\beta \in \ss_4} \frac{1}{2} \#((\R\R_2(d_1,d_2,r),\wt{\d}),(\beta(h),\beta(k)),(\beta(p),\beta(q)))
\big ]
\end{eqnarray*}
\end{theorem} 
\bpf We consider the moduli space $\mbar_{1,4}(r,d)^*$, where there is a one-to-one correspondence 
$\mu : \text{set of marked points} \to \{h,k,p,q\}  $. Let $\F \subset \mbar_{1,4}(r,d)^*$ be the
subfamily cut out by the constraint $\wt{\d}$ and the condition that 
$\gamma(P_i) \in \mu(P_i) $ for each mark point $P_i$. $\F$ is either
empty or $2$-dimensional. If $\F$ is empty then all the summands in the equation
are $0$ so there is nothing to prove. Assume $\F$ is $2$-dimensional. Let $\G$
be the pushforward of $\F$ under the forgetful morphism $\mbar_{1,4}(r,d) \to \mof$.
We intersect $\G$ with Getzler's relation.\\
- $\G \cap \delta_{2,2}$: If the elliptic component is contracted, then by Proposition
$4.2$, then there are two possibilities for the image curves. It can be
a two-component rational curves intersecting at a tacnode. By simple dimension
counting (it lies in a family of one fewer dimension than that of a rational
cuspidal curve) we see that the image curves (of those in $\F$ that contract
the elliptic component) move in a family of $3$ dimension fewer than $\F$. But 
$\F$ is $2$-dimensional so, this loci must be empty. The other possibility is 
again two-component rational curve, but one has a cusp at the intersection. Dimension
counting again shows that this locus on $\F$ is empty. (Note that if $\F$ has
at least $3$ dimensional , then these loci are generally not empty and are in fact
a divisor on $\F$. However, a contributor to the dimensions of 
these loci comes from varying the two nodes
on the elliptic components, which is irrelevant for enumerative reason). If both
rational components are contracted, then we either get $(2+ \one_{(i=2)})\#(\E(r,d),\d)$
depending on if $p$ is $1$-dimensional. If $p$ is not $1$-dimensional, then we also get
$\#(\E(r,d),\d_1)$. Next, if $1$ of the rational tails is contracted, then we get
$$\frac{1}{4} \sum_{d_1 + d_e = d} \sum_{\beta \in \ss_4} \#((\E\R(d_e,d_1,r),\wt{\d}), (\beta(h) \cap \beta(k)), (\beta(p), \beta(q))$$
The factor ${1}/{4}$ comes from the fact that permuting
marked points on the same rational tail gives us the same stratum. If
non of the rational tails are contracted, then the contribution is
$$\frac{1}{8} \sum_{d_1+d_2+d_e = d} \sum_{\beta \in \ss_4} \#((\R\E\R(d_1,d_e,d_r),\wt{\d}),(\beta(h),\beta(k)),\emptyset, (\beta(p),\beta(q)))$$
The factor $1/8$ comes from the fact that 
permuting marked points on each tail and permuting
the two tails give us the same stratum (given that
we consider all possible degree distributions, contrary
to the situation in the previous remark, where
one of the tail is contracted while the other is not). $12$
is the coefficient of $\delta_{2,2}$ in Getzler's relation. \\
- $\G \cap \delta_{2,3}$, $\G \cap \delta_{2,4}, \G \cap \delta_{3,4}$ can be found
similarly, and these summands explain all the terms until the first
$\n\R$ term. \\
- $\G \cap \delta_{0,3}$:
 If the nodal elliptic component is contracted, then by Proposition
$4.2$, we must have a cuspidal curve joined with a smooth rational
curve by a node. But this locus on $\F$ must be empty since
a smooth rational cuspidal curve is moving in a $2$ dimensions
fewer family than that of a smooth elliptic curve (a general
member of $\F$). The rational tail can not be contracted also :
if it were contracted , then moving the marked points on the rational
tail (one-dimensional family) change the map but does not change the
image curves, meaning than the image curves can not satisfy all
the enumerative constraint. Thus non of the components can be contracted,
and we have the contribution from this stratum is
$$\frac{1}{6} \sum_{d_1+d_n =0} \sum_{\beta \in \ss_4} \frac{1}{2}\#((\n\R(d_n,d_1),\wt{\d}),(\beta(h)),(\beta(k),\beta(p),\beta(q)))$$
The factor $1/6$ comes from the fact that we can permute the marked points
on the rational component, and the factor $1/2$ comes from the fact that
we can permute the branches over the node of the image for maps
in $\n(r,d)$. \\
- $\G \cap \delta_{0,4}$: Argue similarly, the rational tail
can not be contracted. If the nodal elliptic component is not
contracted, then the contribution is 
$$\sum_{d_1+d_n =d} \frac{1}{2}\#((\n\R(d_n,d_1,r),\wt{\d}),\emptyset, (h,k,p,q))$$
If the nodal elliptic component is contracted, then the contribution is
$$\frac{1}{2} \#(\s,\d_c)$$
The factor $1/2$ comes from the fact that any map that contracts
the nodal elliptic component has automorphism group of order $2$. \\
- $\G \cap \delta_{b}$ : Argue as above, we can show that
non of the rational bridges can be contracted. The contribution
is 
$$\frac{1}{8}  \sum_{d_1+d_2 =d}\sum_{\beta \in \ss_4} \frac{1}{2} \#((\R\R_2(d_1,d_2,r),\wt{\d}),(\beta(h),\beta(k)),(\beta(p),\beta(q)))$$
We can permute the marked points on each component and permute 
the two component themshelves, hence the factor $1/8$. The factors
$1/2$ comes from the fact that we can permute two branches over
one of the node of the image for maps in $\R\R_2(r,d_1,d_2).$
\epf 

Except the terms involving $\E(r,d)$, all other terms
can be computed either from algorithms in
\cite{dn} and \cite{dn2} or recursively. For example
$\E\R$ terms $\E\R\R$ and $\R\E\R$ can be computed
using the ``splitting the diagonal" method described
in Section $3$ of \cite{dn}. The term $\E(r,d_i)$ for
$i = 0, 1,2$ can be assumed known by induction
since $\d_i$ is of lower rank than that of $\d$.

 Note that we can choose $p',q'$ being any two
subspaces in $\d$ having same codimensions: they
could be two points, two lines, two planes etc. For
each choice we have a different recursion. Hence a
good check of the formulas is that different choices
giving same numbers. We have confirmed this is true
for all elliptic space curves in $\pba, \pbon, \pnam$ of
degree at most $5$, and so far no contradiction has
been found.
\section{Characteristic numbers of elliptic space curves}
In this section we give a recursive formula counting
elliptic space curves with tangency conditions. First we need
the following lemma.
\begin{lemma}
  Let $\d$ be a constraint, and $\F$ a family of maps
in $\mbar_{1,0}(r,d)^*$ such that $\F = (\mbar_{1,0}(r,d)^*,\d)$
is one-dimensional. Then the maps in $\F$ do not
contract an elliptic component. Moreover, the number of
maps in $\F$ whose images are tangent to a general hyperplane
is given by
$$ I + \frac{dJ}{12} + \sum_{i=0}^{r-1}iR_i$$
where $I$ is the number of maps in $\F$ satisfying
top incident condition, $J$ is the number of maps in  $\F$ whose
domain is a smooth elliptic curve having a fixed $j$-invariant
($j$ is not $0$ or $1728)$, and $R_i$ is the number
of rational tails of domains of maps in $\F$ that
are mapped with degree $i$.
\end{lemma}
\bpf The loci of maps that contract an elliptic component
are at most divisors on $\mbar_{1,0}(r,d)$. Since
we can change such maps by moving the nodes on the elliptic
component without changing the images, the images must
move in a family of at most $2$ dimensions fewer than
those of maps of smooth elliptic curves. Since $\F$
is one-dimensional, and enumerative constraints impose
conditions on the image curves only, this shows those loci must be empty.
Let $\T, \H$ be the divisor on $\F$ corresponding to tangency
condition and top incidence condition respectively. Assume the
$\F$ is represented by the total family $\pi : \C \to \F$ and
a map $\mu: \C \to \pr$. Let $H$ be a general hyperplane in $\pr$.
Then $\D= \mu^{(-1)}H$ is a smooth curve in $\C$ and $\T$ is given
by the branch divisor of the covering $\pi : \D \to \F$. Thus
\begin{eqnarray*}
\T\cdot \F = \pi_*(K_{\D} - \pi^*K_{\F}) = \pi_*(K_{\D} - (K_{\C} + \omega_\pi)_{|\D}) = \pi_*(\D(\D+ \omega_{\pi}))
\end{eqnarray*}
where the last equality follows from adjunction. It is
easy to see that $I = \H \cdot \F = \pi_*(\D^2).$
Let $R$ be the divisor of $\C$ corresponding to
rational tails, and let $\delta_0$ be the divisor
of $\C$ corresponding to nodal elliptics. Then
we have
$$\omega_{\pi} = \frac{\delta_0}{12} + R$$ 
This follows from Theorem $12.1$ in \cite{bpv}, corrected by the rational tails. Since
$\delta_0$ is equivalent to the locus of fibers having fixed $j$-invariant,
we have
$$\T \cdot \F = \pi_*(\D^2) + \pi_*(\D \omega_{\pi}) = I + \frac{dJ}{12} + \sum_{i=1}^{r-1} iR_i$$
Using the lemma, we can easily deduce the following result.
\begin{theorem}
 Let $\d$ be a constraint with $\d(0) >0$. Let $\d'$ be the constraint derived from $\d$ by
removing a tangency condition and adding a top incidence condition. Let $\d''$
be the constraint derived from $\d$ by removing a tangency condition. Then we
have the following equality, provided the left-hand side is finite.
$$\#(\E(r,d),\d) = \#(\E(d,r),\d') + \frac{d}{12}\#(\J(r,d),\d') + \sum_{d_1+d_e =d}d_1\#(\E\R(d_e,d_1),\d'').$$
\end{theorem}
 This gives a recursive formula for all characteristic numbers 
of elliptic space curves with at least one tangency
condition. The characteristic numbers $\#(\J(r,d),\d')$ were
computed in \cite{dn}.

\section{Numerical examples}
 In this section we give numerical
examples of characteristic numbers of
elliptic curves of low degree in
$\mathbb P^2, \pba, \pbon, \pnam.$ For
elliptic curves in $\mathbb P^2$ and $\pba$
of low degrees, we give all characteristic
numbers. For curves in higher dimensional
projective spaces, we give a random sample
of these numbers as there are too many possible
characteristic numbers (in addition to
long running time). Note that degree $2$
elliptic curves are understood as degree $2$
covers of $\mathbb P^1$.

$$ $$
$$ $$
$$ $$

\begin{center} 
 \begin{tabular} {| r |  r | r | r | r | r |} 
 \hline
 $\d$ & $E_{\d}$  &$ \d$ & $ E_{\d}$  &$ \d$ & $ E_{\d}$ \\ \hline
$( 0, 6  )$ &0& $( 1, 5  )$ &0& $( 2, 4  )$ &0 \\ 
 \hline
$( 3, 3  )$ &0& $( 4, 2  )$ &2& $( 5, 1  )$ &10 \\ 
 \hline
$( 6, 0  )$ &45/2& & & &  \\ 
 \hline
 \end{tabular} 
 \end{center} 
 \begin{center}  
 Table $1.$ Degree $2$ elliptic plane curves.
 \end{center}

 \begin{center} 
 \begin{tabular} {| r |  r | r | r | r | r |} 
 \hline
 $\d$ & $E_{\d}$  &$ \d$ & $ E_{\d}$  &$ \d$ & $ E_{\d}$ \\ \hline
$( 0, 9 )$ &1& $( 1, 8 )$ &4& $( 2, 7 )$ &16 \\ 
 \hline
$( 3, 6 )$ &64& $( 4, 5 )$ &256& $( 5, 4 )$ &976 \\ 
 \hline
$( 6, 3 )$ &3424& $( 7, 2 )$ &9766& $( 8, 1 )$ &21004 \\ 
 \hline
$( 9, 0 )$ &33616& & & &  \\ 
 \hline
 \end{tabular} 
 \end{center} 
 \begin{center}  
 Table $2.$ Degree $3$ elliptic plane curves.
 \end{center}

\begin{center} 
 \begin{tabular} {| r |  r | r | r | r | r |} 
 \hline
 $\d$ & $E_{\d}$  &$ \d$ & $ E_{\d}$  &$ \d$ & $ E_{\d}$ \\ \hline
$( 0, 12 )$ &225& $( 1, 11 )$ &1010& $( 2, 10 )$ &4396 \\ 
 \hline
$( 3, 9 )$ &18432& $( 4, 8 )$ &73920& $( 5, 7 )$ &280560 \\ 
 \hline
$( 6, 6 )$ &994320& $( 7, 5 )$ &3230956& $( 8, 4 )$ &9409052 \\ 
 \hline
$( 9, 3 )$ &23771160& $( 10, 2 )$ &50569520& $( 11, 1 )$ &89120080 \\ 
 \hline
$( 12, 0 )$ &129996216& & & &  \\ 
 \hline
 \end{tabular} 
 \end{center} 
 \begin{center}  
 Table $3.$  Degree $4$ elliptic plane curves.
 \end{center}

\begin{center} 
 \begin{tabular} {| r |  r | r | r | r | r |} 
 \hline
 $\d$ & $E_{\d}$  &$ \d$ & $ E_{\d}$  &$ \d$ & $ E_{\d}$ \\ \hline
$( 0, 15 )$ &87192& $( 1, 14 )$ &411376& $( 2, 13 )$ &1873388 \\ 
 \hline
$( 3, 12 )$ &8197344& $( 4, 11 )$ &34294992& $( 5, 10 )$ &136396752 \\ 
 \hline
$( 6, 9 )$ &512271756& $( 7, 8 )$ &1802742368& $( 8, 7 )$ &5889847264 \\ 
 \hline
$( 9, 6 )$ &17668868832& $( 10, 5 )$ &48034104112& $( 11, 4 )$ &116575540736 \\ 
 \hline
$( 12, 3 )$ &248984451648& $( 13, 2 )$ &463227482784& $( 14, 1 )$ &747546215472 \\ 
 \hline
$( 15, 0 )$ &1048687299072& & & &  \\ 
 \hline
 \end{tabular} 
 \end{center} 
 \begin{center}  
 Table $4.$ Degree $5$ elliptic plane curves.
 \end{center}

\begin{center} 
 \begin{tabular} {| r |  r | r | r | r | r |} 
 \hline
 $\d$ & $E_{\d}$  &$ \d$ & $ E_{\d}$  &$ \d$ & $ E_{\d}$ \\ \hline
$( 0, 18 )$ &57435240& $( 1, 17 )$ &278443920& $( 2, 16 )$ &1304259360 \\ 
 \hline
$( 3, 15 )$ &5884715280& $( 4, 14 )$ &25491474432& $( 5, 13 )$ &105633321120 \\ 
 \hline
$( 6, 12 )$ &417060737040& $( 7, 11 )$ &1561852784760& $( 8, 10 )$ &5519825870880 \\ 
 \hline
$( 9, 9 )$ &18304445284032& $( 10, 8 )$ &56582281200000& $( 11, 7 )$ &161827650576960 \\ 
 \hline
$( 12, 6 )$ &424685965762560& $( 13, 5 )$ &1013734555246080& $( 14, 4 )$ &2182871531466432 \\ 
 \hline
$( 15, 3 )$ &4212351284630880& $( 16, 2 )$ &7256549594647680& $( 17, 1 )$ &11151931379093760 \\ 
 \hline
$( 18, 0 )$ &15327503832362880& & & &  \\ 
 \hline
 \end{tabular} 
 \end{center} 
 \begin{center}  
 Table $5.$ Degree $6$ elliptic plane curves. 
 \end{center}

\begin{center} 
 \begin{tabular} {| r |  r | r | r | r | r |} 
 \hline
 $\d$ & $E_{\d}$  &$ \d$ & $ E_{\d}$  &$ \d$ & $ E_{\d}$ \\ \hline
$( 0, 0, 6 )$ &0& $( 0, 2, 5 )$ &0& $( 0, 4, 4 )$ &0 \\ 
 \hline
$( 0, 6, 3 )$ &1& $( 0, 8, 2 )$ &14& $( 0, 10, 1 )$ &150 \\ 
 \hline
$( 0, 12, 0 )$ &1500& $( 1, 1, 5 )$ &0& $( 1, 3, 4 )$ &0 \\ 
 \hline
$( 1, 5, 3 )$ &4& $( 1, 7, 2 )$ &50& $( 1, 9, 1 )$ &498 \\ 
 \hline
$( 1, 11, 0 )$ &4740& $( 2, 0, 5 )$ &0& $( 2, 2, 4 )$ &0 \\ 
 \hline
$( 2, 4, 3 )$ &16& $( 2, 6, 2 )$ &176& $( 2, 8, 1 )$ &1620 \\ 
 \hline
$( 2, 10, 0 )$ &14640& $( 3, 1, 4 )$ &0& $( 3, 3, 3 )$ &64 \\ 
 \hline
$( 3, 5, 2 )$ &608& $( 3, 7, 1 )$ &5136& $( 3, 9, 0 )$ &43944 \\ 
 \hline
$( 4, 0, 4 )$ &0& $( 4, 2, 3 )$ &256& $( 4, 4, 2 )$ &2048 \\ 
 \hline
$( 4, 6, 1 )$ &15744& $( 4, 8, 0 )$ &127104& $( 5, 1, 3 )$ &976 \\ 
 \hline
$( 5, 3, 2 )$ &6464& $( 5, 5, 1 )$ &45040& $( 5, 7, 0 )$ &342720 \\ 
 \hline
$( 6, 0, 3 )$ &3424& $( 6, 2, 2 )$ &18560& $( 6, 4, 1 )$ &116768 \\ 
 \hline
$( 6, 6, 0 )$ &836480& $( 7, 1, 2 )$ &47936& $( 7, 3, 1 )$ &269440 \\ 
 \hline
$( 7, 5, 0 )$ &1809040& $( 8, 0, 2 )$ &114248& $( 8, 2, 1 )$ &553176 \\ 
 \hline
$( 8, 4, 0 )$ &3439024& $( 9, 1, 1 )$ &1024404& $( 9, 3, 0 )$ &5768584 \\ 
 \hline
$( 10, 0, 1 )$ &1774680& $( 10, 2, 0 )$ &8656240& $( 11, 1, 0 )$ &11875120 \\ 
 \hline
$( 12, 0, 0 )$ &15480640& & & &  \\ 
 \hline
 \end{tabular} 
 \end{center} 
 \begin{center}  
 Table $6.$ Degree $3$ elliptic curves in $\pba$.
 \end{center}
\begin{center}

 \begin{tabular} {| r |  r | r | r | r | r |} 
 \hline
 $\d$ & $E_{\d}$  &$ \d$ & $ E_{\d}$  &$ \d$ & $ E_{\d}$ \\ \hline
$( 0, 0, 8 )$ &1& $( 0, 2, 7 )$ &4& $( 0, 4, 6 )$ &32 \\ 
 \hline
$( 0, 6, 5 )$ &310& $( 0, 8, 4 )$ &3220& $( 0, 10, 3 )$ &34674 \\ 
 \hline
$( 0, 12, 2 )$ &385656& $( 0, 14, 1 )$ &4436268& $( 0, 16, 0 )$ &52832040 \\ 
 \hline
$( 1, 1, 7 )$ &12& $( 1, 3, 6 )$ &96& $( 1, 5, 5 )$ &920 \\ 
 \hline
$( 1, 7, 4 )$ &9408& $( 1, 9, 3 )$ &99270& $( 1, 11, 2 )$ &1081968 \\ 
 \hline
$( 1, 13, 1 )$ &12224484& $( 1, 15, 0 )$ &143419320& $( 2, 0, 7 )$ &36 \\ 
 \hline
$( 2, 2, 6 )$ &288& $( 2, 4, 5 )$ &2720& $( 2, 6, 4 )$ &27312 \\ 
 \hline
$( 2, 8, 3 )$ &281004& $( 2, 10, 2 )$ &2988144& $( 2, 12, 1 )$ &33049512 \\ 
 \hline
$( 2, 14, 0 )$ &381061200& $( 3, 1, 6 )$ &864& $( 3, 3, 5 )$ &8000 \\ 
 \hline
$( 3, 5, 4 )$ &78656& $( 3, 7, 3 )$ &783840& $( 3, 9, 2 )$ &8087616 \\ 
 \hline
$( 3, 11, 1 )$ &87221808& $( 3, 13, 0 )$ &985671936& $( 4, 0, 6 )$ &2592 \\ 
 \hline
$( 4, 2, 5 )$ &23360& $( 4, 4, 4 )$ &224256& $( 4, 6, 3 )$ &2145024 \\ 
 \hline
$( 4, 8, 2 )$ &21331136& $( 4, 10, 1 )$ &223311840& $( 4, 12, 0 )$ &2466111936 \\ 
 \hline
$( 5, 1, 5 )$ &67440& $( 5, 3, 4 )$ &630720& $( 5, 5, 3 )$ &5721424 \\ 
 \hline
$( 5, 7, 2 )$ &54410016& $( 5, 9, 1 )$ &550239168& $( 5, 11, 0 )$ &5919868800 \\ 
 \hline
$( 6, 0, 5 )$ &191760& $( 6, 2, 4 )$ &1743488& $( 6, 4, 3 )$ &14766080 \\ 
 \hline
$( 6, 6, 2 )$ &133095808& $( 6, 8, 1 )$ &1293435904& $( 6, 10, 0 )$ &13514355840 \\ 
 \hline
$( 7, 1, 4 )$ &4724272& $( 7, 3, 3 )$ &36626544& $( 7, 5, 2 )$ &309751664 \\ 
 \hline
$( 7, 7, 1 )$ &2876272592& $( 7, 9, 0 )$ &29088348480& $( 8, 0, 4 )$ &12532016 \\ 
 \hline
$( 8, 2, 3 )$ &86940920& $( 8, 4, 2 )$ &681603936& $( 8, 6, 1 )$ &6007997008 \\ 
 \hline
$( 8, 8, 0 )$ &58587710176& $( 9, 1, 3 )$ &197671204& $( 9, 3, 2 )$ &1413728496 \\ 
 \hline
$( 9, 5, 1 )$ &11731399560& $( 9, 7, 0 )$ &109792714096& $( 10, 0, 3 )$ &435015624 \\ 
 \hline
$( 10, 2, 2 )$ &2767555376& $( 10, 4, 1 )$ &21376596768& $( 10, 6, 0 )$ &190821802560 \\ 
 \hline
$( 11, 1, 2 )$ &5150502848& $( 11, 3, 1 )$ &36418237824& $( 11, 5, 0 )$ &307505812160 \\ 
 \hline
$( 12, 0, 2 )$ &9269345984& $( 12, 2, 1 )$ &58355286272& $( 12, 4, 0 )$ &460737967360 \\ 
 \hline
$( 13, 1, 1 )$ &88904673408& $( 13, 3, 0 )$ &645526598016& $( 14, 0, 1 )$ &131356680480 \\ 
 \hline
$( 14, 2, 0 )$ &853096310656& $( 15, 1, 0 )$ &1076343432320& $( 16, 0, 0 )$ &1321684733280 \\ 
 \hline
 \end{tabular} 
 \end{center} 
 \begin{center}  
 Table $7.$ Degree $4$ elliptic curves in $\pba$.
 \end{center}

\begin{center} 
 \begin{tabular} {| r |  r | r | r | r | r |} 
 \hline
 $\d$ & $E_{\d}$  &$ \d$ & $ E_{\d}$  &$ \d$ & $ E_{\d}$ \\ \hline
$( 0, 0, 10 )$ &42& $( 0, 2, 9 )$ &354& $( 0, 4, 8 )$ &3492 \\ 
 \hline
$( 0, 6, 7 )$ &38049& $( 0, 8, 6 )$ &441654& $( 0, 10, 5 )$ &5378454 \\ 
 \hline
$( 0, 12, 4 )$ &68292324& $( 0, 14, 3 )$ &901654884& $( 0, 16, 2 )$ &12358163808 \\ 
 \hline
$( 0, 18, 1 )$ &175599635328& $( 0, 20, 0 )$ &2583319387968& $( 1, 1, 9 )$ &1094 \\ 
 \hline
$( 1, 3, 8 )$ &10476& $( 1, 5, 7 )$ &111774& $( 1, 7, 6 )$ &1271974 \\ 
 \hline
$( 1, 9, 5 )$ &15194034& $( 1, 11, 4 )$ &189441324& $( 1, 13, 3 )$ &2460171444 \\ 
 \hline
$( 1, 15, 2 )$ &33232962528& $( 1, 17, 1 )$ &466363099008& $( 1, 19, 0 )$ &6789367904448 \\ 
 \hline
$( 2, 0, 9 )$ &3340& $( 2, 2, 8 )$ &31120& $( 2, 4, 7 )$ &324980 \\ 
 \hline
$( 2, 6, 6 )$ &3618784& $( 2, 8, 5 )$ &42296604& $( 2, 10, 4 )$ &516526368 \\ 
 \hline
$( 2, 12, 3 )$ &6582996888& $( 2, 14, 2 )$ &87478828368& $( 2, 16, 1 )$ &1210578237888 \\ 
 \hline
$( 2, 18, 0 )$ &17419712523648& $( 3, 1, 8 )$ &91560& $( 3, 3, 7 )$ &935136 \\ 
 \hline
$( 3, 5, 6 )$ &10164264& $( 3, 7, 5 )$ &115886944& $( 3, 9, 4 )$ &1381799016 \\ 
 \hline
$( 3, 11, 3 )$ &17235919176& $( 3, 13, 2 )$ &224820084336& $( 3, 15, 1 )$ &3062777447088 \\ 
 \hline
$( 3, 17, 0 )$ &43504838611968& $( 4, 0, 8 )$ &267008& $( 4, 2, 7 )$ &2663824 \\ 
 \hline
$( 4, 4, 6 )$ &28172256& $( 4, 6, 5 )$ &312141824& $( 4, 8, 4 )$ &3619891072 \\ 
 \hline
$( 4, 10, 3 )$ &44047594080& $( 4, 12, 2 )$ &562512674880& $( 4, 14, 1 )$ &7529361687168 \\ 
 \hline
$( 4, 16, 0 )$ &105419849015808& $( 5, 1, 7 )$ &7515344& $( 5, 3, 6 )$ &77029344 \\ 
 \hline
$( 5, 5, 5 )$ &825583024& $( 5, 7, 4 )$ &9266866944& $( 5, 9, 3 )$ &109570881504 \\ 
 \hline
$( 5, 11, 2 )$ &1365937000128& $( 5, 13, 1 )$ &17924770819968& $( 5, 15, 0 )$ &246982965815808 \\ 
 \hline
$( 6, 0, 7 )$ &21015744& $( 6, 2, 6 )$ &207770304& $( 6, 4, 5 )$ &2142245344 \\ 
 \hline
$( 6, 6, 4 )$ &23132708544& $( 6, 8, 3 )$ &264540003744& $( 6, 10, 2 )$ &3208327374528 \\ 
 \hline
$( 6, 12, 1 )$ &41175810566208& $( 6, 14, 0 )$ &557379334146048& $( 7, 1, 6 )$ &553229344 \\ 
 \hline
$( 7, 3, 5 )$ &5451732864& $( 7, 5, 4 )$ &56197250608& $( 7, 7, 3 )$ &618092874928 \\ 
 \hline
$( 7, 9, 2 )$ &7264503650688& $( 7, 11, 1 )$ &90935869598208& $( 7, 13, 0 )$ &1207034339515008 \\ 
 \hline
$( 8, 0, 6 )$ &1456801224& $( 8, 2, 5 )$ &13615688984& $( 8, 4, 4 )$ &132655331216 \\ 
 \hline
$( 8, 6, 3 )$ &1393893072176& $( 8, 8, 2 )$ &15805591490496& $( 8, 10, 1 )$ &192394855764288 \\ 
 \hline
$( 8, 12, 0 )$ &2498937360190848& $( 9, 1, 5 )$ &33442096324& $( 9, 3, 4 )$ &304042170552 \\ 
 \hline
$( 9, 5, 3 )$ &3027600730440& $( 9, 7, 2 )$ &32950822741600& $( 9, 9, 1 )$ &388703233243008 \\ 
 \hline
$( 9, 11, 0 )$ &4928984795256768& $( 10, 0, 5 )$ &81107025144& $( 10, 2, 4 )$ &677080161264 \\ 
 \hline
$( 10, 4, 3 )$ &6326216895824& $( 10, 6, 2 )$ &65682090217248& $( 10, 8, 1 )$ &747920741035008 \\ 
 \hline
$( 10, 10, 0 )$ &9234940120602048& $( 11, 1, 4 )$ &1469692262864& $( 11, 3, 3 )$ &12719562619344 \\ 
 \hline
$( 11, 5, 2 )$ &125046496125728& $( 11, 7, 1 )$ &1368071583917408& $( 11, 9, 0 )$ &16398336070362048 \\ 
 \hline
$( 12, 0, 4 )$ &3132954358848& $( 12, 2, 3 )$ &24665727975168& $( 12, 4, 2 )$ &227432721556608 \\ 
 \hline
$( 12, 6, 1 )$ &2377092153260928& $( 12, 8, 0 )$ &27559370814663168& $( 13, 1, 3 )$ &46380697328576 \\ 
 \hline
$( 13, 3, 2 )$ &395961856525056& $( 13, 5, 1 )$ &3925472901986688& $( 13, 7, 0 )$ &43826076405481728 \\ 
 \hline
$( 14, 0, 3 )$ &85453000115200& $( 14, 2, 2 )$ &662634098949120& $( 14, 4, 1 )$ &6173065500610048 \\ 
 \hline
$( 14, 6, 0 )$ &66008781304150528& $( 15, 1, 2 )$ &1073447856471168& $( 15, 3, 1 )$ &9277390256888448 \\ 
 \hline
$( 15, 5, 0 )$ &94369223149298688& $( 16, 0, 2 )$ &1703768515184128& $( 16, 2, 1 )$ &13397726490436608 \\ 
 \hline
$( 16, 4, 0 )$ &128504536404712448& $( 17, 1, 1 )$ &18739152863106048& $( 17, 3, 0 )$ &167464728764784128 \\ 
 \hline
$( 18, 0, 1 )$ &25701019852524288& $( 18, 2, 0 )$ &210144272908228608& $( 19, 1, 0 )$ &255958547477177088 \\ 
 \hline
$( 20, 0, 0 )$ &306095919912649728& & & &  \\ 
 \hline
 \end{tabular} 
 \end{center} 
 \begin{center}  
 Table $8.$ Degree $5$ elliptic curves in $\pba$. 
 \end{center}

\begin{center} 
 \begin{tabular} {| r |  r | r | r |} 
 \hline
 $\d$ & $E_{\d}$  &$ \d$ & $ E_{\d}$ \\ \hline
$( 0, 0, 3, 3 )$ &0& $( 0, 2, 2, 3 )$ &0 \\ 
 \hline
$( 0, 3, 3, 2 )$ &0& $( 0, 6, 0, 3 )$ &1 \\ 
 \hline
$( 0, 7, 1, 2 )$ &14& $( 0, 8, 2, 1 )$ &222 \\ 
 \hline
$( 0, 9, 0, 2 )$ &114& $( 1, 4, 2, 2 )$ &4 \\ 
 \hline
$( 1, 4, 5, 0 )$ &190& $( 1, 6, 4, 0 )$ &1488 \\ 
 \hline
$( 2, 0, 2, 3 )$ &0& $( 2, 0, 5, 1 )$ &0 \\ 
 \hline
$( 2, 1, 3, 2 )$ &0& $( 2, 3, 2, 2 )$ &16 \\ 
 \hline
$( 2, 3, 5, 0 )$ &640& $( 2, 6, 2, 1 )$ &2280 \\ 
 \hline
$( 2, 7, 3, 0 )$ &31044& $( 3, 4, 4, 0 )$ &13680 \\ 
 \hline
$( 3, 8, 2, 0 )$ &536304& $( 4, 0, 4, 1 )$ &512 \\ 
 \hline
$( 4, 2, 3, 1 )$ &3328& $( 4, 3, 1, 2 )$ &2048 \\ 
 \hline
$( 4, 9, 1, 0 )$ &8041776& $( 5, 2, 1, 2 )$ &6464 \\ 
 \hline
$( 5, 2, 4, 0 )$ &101408& $( 5, 4, 0, 2 )$ &31624 \\ 
 \hline
$( 5, 4, 3, 0 )$ &582888& $( 6, 0, 3, 1 )$ &26848 \\ 
 \hline
$( 6, 2, 2, 1 )$ &144400& $( 6, 3, 0, 2 )$ &83312 \\ 
 \hline
$( 6, 7, 1, 0 )$ &39311360& $( 7, 0, 4, 0 )$ &497216 \\ 
 \hline
$( 7, 2, 0, 2 )$ &203968& $( 7, 6, 1, 0 )$ &76501840 \\ 
 \hline
$( 9, 0, 3, 0 )$ &8583182& $( 9, 4, 1, 0 )$ &224882706 \\ 
 \hline
$( 9, 6, 0, 0 )$ &1130248810& $( 10, 0, 1, 1 )$ &10539980 \\ 
 \hline
$( 11, 1, 0, 1 )$ &72275990& $( 15, 0, 0, 0 )$ &5552993600 \\ 
 \hline
 \end{tabular} 
 \end{center} 
 \begin{center}  
 Table $9.$ Some characteristic numbers of elliptic space curves of degree $3$ in $\pbon.$
 \end{center}

\begin{center} 
 \begin{tabular} {| r |  r | r | r |} 
 \hline
 $\d$ & $E_{\d}$  &$ \d$ & $ E_{\d}$ \\ \hline
$( 0, 0, 7, 2 )$ &29& $( 0, 1, 2, 5 )$ &0 \\ 
 \hline
$( 0, 2, 3, 4 )$ &4& $( 0, 10, 2, 2 )$ &1004916 \\ 
 \hline
$( 0, 13, 2, 1 )$ &144007483& $( 1, 0, 5, 3 )$ &38 \\ 
 \hline
$( 1, 3, 5, 2 )$ &4860& $( 1, 4, 6, 1 )$ &81196 \\ 
 \hline
$( 1, 5, 1, 4 )$ &920& $( 1, 7, 0, 4 )$ &9408 \\ 
 \hline
$( 1, 7, 6, 0 )$ &11043310& $( 1, 8, 1, 3 )$ &156968 \\ 
 \hline
$( 1, 8, 4, 1 )$ &5342984& $( 1, 10, 3, 1 )$ &43094568 \\ 
 \hline
$( 2, 7, 1, 3 )$ &420000& $( 2, 7, 4, 1 )$ &13098172 \\ 
 \hline
$( 2, 12, 3, 0 )$ &13916950104& $( 2, 13, 1, 1 )$ &6377111884 \\ 
 \hline
$( 2, 14, 2, 0 )$ &113231319632& $( 3, 0, 1, 5 )$ &0 \\ 
 \hline
$( 3, 1, 2, 4 )$ &864& $( 3, 1, 8, 0 )$ &958400 \\ 
 \hline
$( 3, 3, 7, 0 )$ &7666224& $( 3, 4, 2, 3 )$ &139312 \\ 
 \hline
$( 4, 1, 3, 3 )$ &47136& $( 4, 3, 5, 1 )$ &9606144 \\ 
 \hline
$( 4, 4, 3, 2 )$ &5203072& $( 5, 1, 1, 4 )$ &67440 \\ 
 \hline
$( 5, 3, 3, 2 )$ &12588560& $( 5, 6, 0, 3 )$ &49299816 \\ 
 \hline
$( 6, 0, 4, 2 )$ &4092688& $( 6, 6, 4, 0 )$ &35732553632 \\ 
 \hline
$( 7, 8, 1, 1 )$ &257686679704& $( 7, 13, 0, 0 )$ &214317637545920 \\ 
 \hline
$( 8, 2, 2, 2 )$ &965893920& $( 9, 2, 3, 1 )$ &19635766224 \\ 
 \hline
$( 11, 7, 1, 0 )$ &194300922530150& $( 12, 3, 1, 1 )$ &3865044261124 \\ 
 \hline
$( 12, 5, 0, 1 )$ &23269075459780& $( 13, 1, 0, 2 )$ &709454579680 \\ 
 \hline
$( 13, 7, 0, 0 )$ &2593390019349960& $( 15, 3, 1, 0 )$ &680573664677760 \\ 
 \hline
$( 17, 1, 1, 0 )$ &1046525920177280& $( 17, 3, 0, 0 )$ &6253213828581120 \\ 
 \hline
 \end{tabular} 
 \end{center} 
 \begin{center}  
 Table $10.$ Some characteristic numbers of elliptic space curves of degree $4$ in $\pbon.$
 \end{center}
\begin{center} 
 \begin{tabular} {| r |  r | r | r |} 
 \hline
 $\d$ & $E_{\d}$  &$ \d$ & $ E_{\d}$ \\ \hline
$( 0, 4, 9, 1 )$ &17812920& $( 1, 0, 0, 8 )$ &8 \\ 
 \hline
$( 1, 4, 7, 2 )$ &22141610& $( 1, 13, 1, 3 )$ &83534718240 \\ 
 \hline
$( 2, 1, 8, 2 )$ &5574120& $( 2, 10, 2, 3 )$ &21135644112 \\ 
 \hline
$( 2, 19, 2, 0 )$ &157704159607499400& $( 3, 6, 2, 4 )$ &317850984 \\ 
 \hline
$( 3, 13, 0, 3 )$ &3852663999504& $( 4, 6, 3, 3 )$ &12466665712 \\ 
 \hline
$( 4, 12, 3, 1 )$ &324662678747600& $( 4, 13, 4, 0 )$ &6438834486545696 \\ 
 \hline
$( 5, 2, 9, 0 )$ &152170558528& $( 5, 8, 3, 2 )$ &4059955629840 \\ 
 \hline
$( 5, 11, 0, 3 )$ &17870726781040& $( 6, 1, 3, 4 )$ &495750960 \\ 
 \hline
$( 6, 7, 6, 0 )$ &267364637315200& $( 6, 10, 0, 3 )$ &37308946152960 \\ 
 \hline
$( 6, 13, 0, 2 )$ &5661536113375616& $( 7, 9, 0, 3 )$ &76189759273400 \\ 
 \hline
$( 7, 10, 1, 2 )$ &1261503043660160& $( 7, 13, 1, 1 )$ &202527736950333560 \\ 
 \hline
$( 8, 8, 0, 3 )$ &152112968567744& $( 8, 14, 0, 1 )$ &3296407750001057840 \\ 
 \hline
$( 10, 4, 1, 3 )$ &75941963996580& $( 11, 2, 0, 4 )$ &13254336373174 \\ 
 \hline
$( 12, 2, 4, 1 )$ &4927564446146120& $( 12, 3, 5, 0 )$ &67069699366126280 \\ 
 \hline
$( 13, 0, 3, 2 )$ &686078803431960& $( 13, 6, 0, 2 )$ &305992051575799496 \\ 
 \hline
$( 13, 9, 0, 1 )$ &36003235532769890160& $( 14, 11, 0, 0 )$ &7621648299681170363680 \\ 
 \hline
$( 15, 7, 0, 1 )$ &76323334247130289600& $( 17, 1, 2, 1 )$ &2598949329007010304 \\ 
 \hline
$( 18, 1, 0, 2 )$ &2695498962098743296& &  \\ 
 \hline
 \end{tabular} 
 \end{center} 
 \begin{center}  
 Table $11.$ Some characteristic numbers of elliptic space curves of degree $5$ in $\pbon.$
 \end{center}

\begin{center} 
 \begin{tabular} {| r |  r | r | r |} 
 \hline
 $\d$ & $E_{\d}$  &$ \d$ & $ E_{\d}$ \\ \hline
$( 0, 0, 6, 2, 0 )$ &0& $( 0, 4, 4, 2, 0 )$ &112 \\ 
 \hline
$( 0, 7, 2, 1, 1 )$ &294& $( 1, 2, 2, 1, 2 )$ &0 \\ 
 \hline
$( 1, 3, 2, 2, 1 )$ &8& $( 1, 3, 5, 0, 1 )$ &190 \\ 
 \hline
$( 1, 5, 3, 2, 0 )$ &2694& $( 2, 0, 1, 2, 2 )$ &0 \\ 
 \hline
$( 2, 0, 6, 0, 1 )$ &80& $( 2, 4, 0, 4, 0 )$ &432 \\ 
 \hline
$( 3, 5, 0, 2, 1 )$ &7008& $( 3, 5, 1, 0, 2 )$ &3696 \\ 
 \hline
$( 4, 0, 5, 0, 1 )$ &6400& $( 4, 2, 0, 0, 3 )$ &256 \\ 
 \hline
$( 5, 0, 0, 3, 1 )$ &976& $( 5, 0, 3, 1, 1 )$ &13536 \\ 
 \hline
$( 5, 1, 1, 2, 1 )$ &10000& $( 5, 5, 0, 0, 2 )$ &135504 \\ 
 \hline
$( 6, 0, 4, 0, 1 )$ &235888& $( 6, 5, 2, 1, 0 )$ &63034080 \\ 
 \hline
$( 6, 7, 1, 1, 0 )$ &307313568& $( 7, 0, 4, 1, 0 )$ &5158008 \\ 
 \hline
$( 8, 10, 0, 0, 0 )$ &160210908480& $( 9, 9, 0, 0, 0 )$ &242570957664 \\ 
 \hline
$( 12, 6, 0, 0, 0 )$ &617311046976& $( 13, 2, 0, 1, 0 )$ &31799786256 \\ 
 \hline
$( 14, 1, 0, 1, 0 )$ &40866566832& $( 14, 2, 1, 0, 0 )$ &219076067112 \\ 
 \hline
 \end{tabular} 
 \end{center} 
 \begin{center}  
 Table $12.$ Some characteristic numbers of elliptic space curves of degree $3$ in $\pnam.$
 \end{center}
\begin{center} 
 \begin{tabular} {| r |  r | r | r |} 
 \hline
 $\d$ & $E_{\d}$  &$ \d$ & $ E_{\d}$ \\ \hline
$( 0, 3, 7, 1, 1 )$ &150696& $( 0, 8, 0, 0, 4 )$ &3220 \\ 
 \hline
$( 0, 9, 0, 1, 3 )$ &57960& $( 0, 10, 3, 0, 2 )$ &18866664 \\ 
 \hline
$( 0, 17, 0, 1, 1 )$ &139530931995& $( 0, 20, 0, 0, 1 )$ &11150743642205 \\ 
 \hline
$( 1, 1, 2, 2, 3 )$ &12& $( 1, 6, 3, 1, 2 )$ &585468 \\ 
 \hline
$( 1, 10, 0, 3, 1 )$ &56249988& $( 2, 0, 2, 6, 0 )$ &7544 \\ 
 \hline
$( 2, 0, 5, 0, 3 )$ &1444& $( 2, 1, 7, 1, 1 )$ &856860 \\ 
 \hline
$( 2, 5, 4, 3, 0 )$ &75690880& $( 2, 14, 4, 0, 0 )$ &36244480608680 \\ 
 \hline
$( 3, 0, 3, 5, 0 )$ &203200& $( 3, 4, 1, 5, 0 )$ &10274688 \\ 
 \hline
$( 3, 4, 2, 3, 1 )$ &6044208& $( 3, 11, 3, 0, 1 )$ &299599503152 \\ 
 \hline
$( 4, 1, 6, 1, 1 )$ &30574272& $( 4, 3, 7, 1, 0 )$ &5402075552 \\ 
 \hline
$( 4, 5, 2, 1, 2 )$ &56706560& $( 5, 10, 3, 1, 0 )$ &22175099750880 \\ 
 \hline
$( 6, 1, 7, 1, 0 )$ &19987349088& $( 6, 9, 1, 1, 1 )$ &1477032997824 \\ 
 \hline
$( 6, 14, 0, 0, 1 )$ &555756395786592& $( 7, 0, 7, 1, 0 )$ &36471958824 \\ 
 \hline
$( 8, 5, 2, 1, 1 )$ &803465851328& $( 9, 1, 4, 2, 0 )$ &479197387872 \\ 
 \hline
$( 9, 11, 0, 0, 1 )$ &2604953865714080& $( 11, 2, 2, 1, 1 )$ &3956581280704 \\ 
 \hline
$( 13, 2, 3, 1, 0 )$ &865530347368728& $( 16, 0, 0, 0, 2 )$ &127384328451776 \\ 
 \hline
 \end{tabular} 
 \end{center} 
 \begin{center}  
 Table $13.$ Some characteristic numbers of elliptic space curves of degree $4$ in $\pnam.$
 \end{center}

\begin{center} 
 \begin{tabular} {| r |  r | r | r |} 
 \hline
 $\d$ & $E_{\d}$  &$ \d$ & $ E_{\d}$ \\ \hline
$( 0, 4, 7, 0, 3 )$ &21370599& $( 0, 5, 0, 3, 4 )$ &13368 \\ 
 \hline
$( 0, 16, 1, 0, 3 )$ &6441662050785& $( 1, 1, 1, 2, 5 )$ &120 \\ 
 \hline
$( 2, 24, 2, 0, 0 )$ &92992972101933954474544& $( 5, 6, 1, 3, 2 )$ &669439268976 \\ 
 \hline
$( 6, 0, 0, 4, 3 )$ &70316304& $( 8, 2, 4, 4, 0 )$ &249642503347264 \\ 
 \hline
$( 8, 9, 0, 3, 1 )$ &59596479340997440& $( 9, 11, 5, 0, 0 )$ &4670806268148474812288 \\ 
 \hline
$( 10, 4, 4, 0, 2 )$ &19278514559525472& $( 11, 2, 0, 3, 2 )$ &275265022278424 \\ 
 \hline
$( 11, 4, 0, 1, 3 )$ &1322240935084886& $( 11, 15, 2, 0, 0 )$ &5731941475630570274830480 \\ 
 \hline
$( 12, 3, 0, 5, 0 )$ &93646029368101232& $( 12, 6, 0, 4, 0 )$ &7985056411618359072 \\ 
 \hline
$( 13, 7, 0, 2, 1 )$ &55781303615787140368& $( 15, 1, 1, 0, 3 )$ &110267644926473616 \\ 
 \hline
 \end{tabular} 
 \end{center} 
 \begin{center}  
 Table $14$ Some characteristic numbers of elliptic space curves of degree $5$ in $\pnam.$
 \end{center}


\begin{thebibliography}{[1]}
\bibitem[AV]{av} D. Avritzer and I. Vainsencher, {\em Compactifying the space of elliptic quartic curves,} in Complex Projective
Geometry, G. Ellingsrud et al. eds., Cambridge U.P.: Cambridge, 1992.
\bibitem[BPV] {bpv} W. Barth, C. Peters, A. Van de Ven, {\em Compact complex surfaces} , Springer-Verlag, New York, 1984.
\bibitem[F]{f} W.Fulton {\em Intersection Theory,} Second Edition, Springer 1996.
\bibitem[FP] {fp} W.Fulton and R. Pandharipande, {\em Notes on 
stable maps and quantum cohomology}, preprint 1996, alg-geom/9608011.
\bibitem[HM] {hm} X. Hernandez, J. M. Miret, {\em The characteristic numbers of cuspidal plane cubics in $\mathbb P^3$} , Bull. Belg.
Math. Soc. Simon Stevin, {\bf 10} (2003) No. 1, 115--124.
\bibitem[HMX] {hmx} X. Hernandez, J. M. Miret and S. Xambo-Descamps, {\em Computing the characteristic numbers of the
variety of nodal plane cubics in $\mathbb P^3$} , J. Symb. Comp. {\bf 42} (2007) 192--202.
\bibitem[I] {eln} E. Ionel, {\em Genus-one enumerative invariants in $\mathbb{P}^n$
 with fixed j-invariant}, Duke Math. J. {\bf 94 (2)} (1998) 279--324.
\bibitem[G] {eg} E. Getzler, {\em Intersection theory on $\mbar_{1,4}$ and elliptic Gromov-Witten invariants},
 J. Amer. Math. Soc. {\bf 10} No. 4 (1997) 973--998.
\bibitem[GKP] {gkp} T. Graber, J. Kock and R. Pandharipande {\em Descendant Invariants and Characteristic Numbers, }
American Journal of Mathematics, Vol. 124, No. {\bf 3}(2002).
\bibitem[KSX] {ksx} S. Kleiman, S. A. Strømme and S. Xambó 
{\em Sketch of a verification of Schubert's number 5819539783680 of twisted cubics,} Lecture Notes in Mathematics, 1987, Volume 1266/1987, 156-180.
\bibitem[N1] {dn} D. Nguyen, {\em Characteristic numbers of elliptic curves with fixed $j-$invariant}, preprint 2011, arXiv:1111.6295. 
\bibitem[N2] {dn2} D. Nguyen, {\em Characteristic numbers of rational cuspidal space curves}, preprint 2011, arXiv:1111.6296.
\bibitem[P1] {idq} R. Pandharipande, {\em Intersection of  $\Q$-divisors on 
Kontsevich's moduli space $\mbar_{0,n}(\mathbb{P}^r,d)$ and enumerative 
geometry}, Trans. Amer. Math. Soc, {\bf 351} (1999), 1481-1505.
\bibitem[P2] {pj} R. Pandharipande, {\em A note on elliptic plane curves with 
fixed $j$-invariant}, Proc. Amer. Math. Soc., {\bf 125}, No. 12, 3471--3479.
\bibitem[P3] {pg}  R. Pandharipande, {\em A geometric construction of Getzler’s elliptic
relation,} Math. Ann. {\bf 313} (1999), 715–729
\bibitem[V1] {ratell} R. Vakil, {\em The enumerative geometry of rational
and elliptic plane curves in projective space}, J. Reine Angew. Math. (Crelle's Journal), {\bf 529} (2000), 101--153.
\bibitem[V2] {char} R. Vakil, {\em  Recursions for characteristic numbers of genus
 one plane curves}, Arkiv for Matematik, {\bf 39} (2001), no. 1, 157--180.
\bibitem[VZ]{vz} R. Vakil, A. Zinger, {\em A desingularization 
of the main component of the moduli space of genus-one stable 
maps to projective space}, Geom. Topol. {\bf 12} (2008), no. 1, 1-95.
\bibitem[Z] {zin} A. Zinger, {\em Enumeration of one-nodal rational curves 
in projective spaces }, Topology {\bf 43} (2004) 793--829.
\end{thebibliography}
\end{document}